\documentclass[12pt]{article}
\usepackage{latexsym}

\def\obr#1 {
\unitlength1mm
\begin{picture}(30,14)(0,0)
#1
\end{picture}
}
\newcommand{\bod}{\circle*{1.5}}

\newcommand{\duk}{\noindent {\bf Proof. }}
\newcommand{\kduk}{\hfill $\Box$\bigskip}
\newcommand{\R}{\mathbf{R}}
\newcommand{\N}{\mathbf{N}}

\newcommand{\ex}{\mathrm{ex}}
\newcommand{\gex}{\mathrm{gex}}

\newtheorem{veta}{Theorem}[section]

\newtheorem{lema}[veta]{Lemma}
\newtheorem{defi}[veta]{Definition}
\newtheorem{prob}[veta]{Problem}

\def\cla#1#2#3#4#5#6{
  {\sc #1, }#2, {\it #3, }{\bf #4 }(#5), #6.}

\def\pre#1#2#3#4#5{
  {\sc #1, }#2, {\it #3, }technical report {\bf #4}, #5.}

\def\kni#1#2#3#4#5{
  {\sc #1, }{\it #2, }#3, #4, #5.}

\def\vsbo#1#2#3#4#5#6#7#8{
  {\sc #1, }#2. In: {#4 (ed.), } {\it #5, } #6, #7,  #8; pp. #3.}

\begin{document}

\author{Martin Klazar\thanks{Department of Applied Mathematics (KAM) and Institute for Theoretical 
Computer Science (ITI), Charles University, Malostransk\'e n\'am\v est\'\i\ 25, 118 00 Praha, 
Czech Republic. ITI is supported by the project LN00A056 of the 
Ministery of Education of the Czech Republic. E-mail: {\tt klazar@kam.mff.cuni.cz}}}
\title{Extremal problems for ordered (hyper)graphs: applications of Davenport--Schinzel sequences}
\date{}

\maketitle
\begin{abstract}
We introduce a containment relation of hypergraphs which respects linear orderings of vertices and investigate 
associated extremal functions. We extend, by means of a more generally applicable theorem, the $n\log n$ 
upper bound on the ordered graph extremal function of $F=(\{1,3\},\{1,5\},\{2,3\},\{2,4\})$ due to F\"uredi to 
the $n(\log n)^2(\log\log n)^3$ upper bound in the hypergraph case. We use Davenport--Schinzel sequences to 
derive almost linear upper bounds in terms of the inverse Ackermann function $\alpha(n)$. We obtain such upper
bounds for the extremal functions of forests consisting of stars whose all centers precede all leaves.  
\end{abstract}

\section{Introduction and motivation}

In this article we shall investigate extremal problems on graphs and hypergraphs of 
the following type. Let $G=([n],E)$ be a simple graph that has the vertex set 
$[n]=\{1,2,\dots,n\}$ and contains no six vertices $1\le v_1<v_2<\cdots<v_6\le n$ such that 
$\{v_1,v_3\},\{v_1,v_5\},\{v_2,v_4\},$ and $\{v_2,v_6\}$ are edges of $G$, that is, $G$ has no {\em ordered\/} 
subgraph of the form
\begin{equation}\label{G0}
G_0=\ \ \obr{\put(0,1){\bod}\put(6,1){\bod}\put(12,1){\bod} 
\put(18,1){\bod}\put(24,1){\bod}\put(30,1){\bod}
\put(6,1){\oval(12,9)[t]}
\put(12,1){\oval(24,13)[t]}
\put(12,1){\oval(12,6)[t]}
\put(18,1){\oval(24,19)[t]}
}
\ \ .
\end{equation}
Determine the maximum possible number $g(n)=|E|$ of edges in $G$.  

What makes this task hard is the linear ordering of $V=[n]$ and the fact that $G_0$ must not appear in $G$ 
only as an ordered subgraph. If we ignore  
the ordering for a while, then the problem asks to determine the maximum number of edges in a simple graph 
$G$ with $n$ vertices and no $2K_{1,2}$ subgraph, and is easily solved. Clearly, if
$G$ has two vertices of degrees $\ge 3$ and $\ge 5$, respectively, or if $G$ has $\ge 6$ vertices of 
degrees $4$ each, then a $2K_{1,2}$ subgraph must appear. Thus if $G$ has a vertex of degree $\ge 5$ and no 
$2K_{1,2}$ subgraph, it has at most $(2(n-1)+n-1)/2=3n/2-1.5$ edges. If all degrees are $\le 4$, the number of 
edges is at most $(3(n-5)+4\cdot 5)/2=3n/2+2.5$. On the other hand, the graph on $[n]$ in which $\deg(n)=n-1$ 
and $[n-1]$ induces a matching with $\lfloor (n-1)/2\rfloor$ edges has $n+\lfloor (n-1)/2\rfloor-1$ edges and no 
$2K_{1,2}$ subgraph. We conclude that in the unordered version of the problem the maximum number of edges 
equals $3n/2+O(1)$.

The ordered version is considerably more difficult. Later in this section we prove that the maximum 
number of edges $g(n)$ satisfies 
\begin{equation}\label{gn}
n\cdot\alpha(n)\ll g(n)\ll n\cdot 2^{(1+o(1))\alpha(n)^2}
\end{equation} 
where $\alpha(n)$ is the inverse Ackermann function. Recall that $\alpha(n)=\min\{m:\ A(m)\ge n\}$ where 
$A(n)=F_n(n)$, the Ackermann function, is defined as follows. We start with $F_1(n)=2n$ and for $i\ge 1$ 
we define $F_{i+1}(n)=F_i(F_i(\ldots F_i(1)\ldots ))$ with $n$ iterations of $F_i$. The function $\alpha(n)$ 
grows to infinity but its growth is extremely slow.
We obtain (\ref{gn}) and some generalizations by 
reductions to {\em Davenport--Schinzel sequences\/} (DS sequences). We continue now with a brief review 
of results on DS sequences that will be needed in the following.  We refer the reader for more information 
and references on DS sequences and their 
applications in computational and combinatorial geometry to Agarwal and Sharir 
\cite{agar_shar}, Klazar \cite{klaz02}, Sharir and Agarwal \cite{shar_agar}, and Valtr \cite{valt99}. 

If $u=a_1a_2\dots a_r$ and $v=b_1b_2\dots b_s$ are two finite sequences (words) over a fixed infinite alphabet
$A$, where $A$ contains $\N=\{1,2,\dots\}$ and also some symbols $a,b,c,d,\dots$, we say that 
$v$ {\em contains\/} $u$ and 
write $v\succ u$ if $v$ has a subsequence $b_{i_1}b_{i_2}\dots b_{i_r}$ such that for every $p$ and $q$ 
we have $a_p=a_q$ if and only if $b_{i_p}=b_{i_q}$. In other words, $v$ has a subsequence that differs from 
$u$ only by an injective renaming of the symbols. For example, $v=ccaaccbaa\succ 22244=u$ because $v$ has 
the subsequence $cccaa$. On the other hand, $ccaaccbaa\not\succ 12121$. A sequence $u=a_1a_2\dots a_r$ is 
called $k$-{\em sparse\/}, where $k\in\N$, if $a_i=a_j,i<j,$ always implies $j-i\ge k$; this means that
every interval in $u$ of length at most $k$ consists of distinct terms. The length $r$ of $u$ is denoted by 
$|u|$. For two integers $a\le b$ we write $[a,b]$ for the interval $\{a,a+1,\dots,b\}$. For two functions 
$f,g:\N\to\R$ the notation $f\ll g$ is synonymous to the $f=O(g)$ notation; it means that $|f(n)|<c|g(n)|$ for 
all $n>n_0$ with a constant $c>0$.

The classical theory of DS sequences 
investigates, for a fixed $s\in\N$, the function $\lambda_s(n)$ that is defined as the maximum length of 
a 2-sparse sequence $v$ over $n$ symbols which does not contain the $s+2$-term alternating sequence 
$ababa\dots$ ($a\ne b$). The notation $\lambda_s(n)$ and the shift $+2$ are due to historical reasons. The
term {\em DS sequences\/} refers to the sequences $v$ not containing a fixed alternating sequence.  
The theory of {\em generalized DS sequences\/} investigates, for a fixed sequence $u$ that uses exactly $k$ 
symbols, the function $\ex(u,n)$ that is defined as the maximum length of a $k$-sparse sequence $v$ such that 
$v$ is over $n$ symbols and $v\not\succ u$. Note that $\ex(u,n)$ extends $\lambda_s(n)$ since 
$\lambda_s(n)=\ex(ababa\dots,n)$ where $ababa\dots$ has length $s+2$. In the definition of 
$\ex(u,n)$ one has to require that $v$ is $k$-sparse because no condition or even 
only $k-1$-sparseness would allow an infinite $v$ with $v\not\succ u$; for example, 
$v=12121212\dots\not\succ abca=u$ and $v$ is 2-sparse (but not 3-sparse). An easy pigeon hole argument shows 
that always $\ex(u,n)<\infty$. 

DS sequences were introduced by Davenport and Schinzel \cite{dave_schi} and strongest bounds on 
$\lambda_s(n)$ for general $s$ were obtained by Agarwal, Sharir and Shor \cite{agar_shar_shor}. 
We need their bound 
\begin{equation}\label{lambda6}
\lambda_6(n)\ll n\cdot 2^{(1+o(1))\alpha(n)^2}
\end{equation}
(recall that $\lambda_6(n)=\ex(abababab,n)$). Hart and Sharir \cite{hart_shar} proved that  
\begin{equation}\label{lambda3}
n\alpha(n)\ll\lambda_3(n)\ll n\alpha(n). 
\end{equation}
In Klazar \cite{klaz92} we proved that if $u$ is a sequence using $k\ge 2$ symbols and $|u|=l\ge 5$, then 
for every $n\in\N$
\begin{equation}\label{mujobec}
\ex(u,n)\le n\cdot k2^{l-3}\cdot (10k)^{2\alpha(n)^{l-4}+8\alpha(n)^{l-5}}.
\end{equation}
It is easy to show that for $k=1$ or $l\le 4$ we have $\ex(u,n)\ll n$. In particular, for the sequence
\begin{equation}\label{ukl}
u(k,l)=12\ldots k12\ldots k\ldots 12\ldots k
\end{equation}
with $l$ segments $12\ldots k$ we have, for every fixed $k\ge 2$ and $l\ge 3$, 
\begin{equation}\label{muj}
\ex(u(k,l),n)\le n\cdot k2^{kl-3}\cdot (10k)^{2\alpha(n)^{kl-4}+8\alpha(n)^{kl-5}}.
\end{equation}
We denote the factor at $n$ in (\ref{muj}) as $\beta(k,l,n)$. Thus
\begin{equation}\label{beta}
\beta(k,l,n)=k2^{kl-3}(10k)^{2\alpha(n)^{kl-4}+8\alpha(n)^{kl-5}}.
\end{equation}

Let us see now how (\ref{lambda6}) and the lower bound in (\ref{lambda3}) imply (\ref{gn}). Let $G=([n],E)$ 
be any simple graph not containing $G_0$ (given in (\ref{G0})) as an ordered subgraph. Consider the sequence 
$$
v=I_1I_2\dots I_n
$$ 
over $[n]$ where $I_i$ is the 
decreasing ordering of the list $\{j:\ \{j,i\}\in E\;\&\;j<i\}$. Note that $I_1=\emptyset$ and $|v|=|E|$. 

\begin{lema}
If $v\succ abababab$ then $G_0$ is an ordered subgraph of $G$.
\end{lema}
\duk
We assume that $v$ has an 8-term alternating subsequence
$$
\dots a_1\dots b_1\dots a_2\dots b_2\dots a_3\dots b_3\dots a_4\dots b_4\dots
$$
where the appearances of two numbers $a\neq b$ are indexed for further discussion. We distinguish two cases. 
If $a<b$ then $a_2$, $b_2$, $a_4$, and $b_4$ lie, respectively, in four distinct intervals $I_p$, $I_q$, 
$I_r$, and $I_s$, $p<q<r<s$, (since every $I_i$ is decreasing) and $b<p$ (since $b_1$ precedes $a_2$). Hence  
$G_0$ is an ordered subgraph of $G$. If $b<a$ then $b_1$, $a_2$, $b_3$, and $a_4$ lie, respectively, in four 
distinct intervals $I_p$, $I_q$, $I_r$, and $I_s$, $p<q<r<s$, and $a<p$. Again, $G_0$ is an ordered subgraph 
of $G$.
\kduk

\noindent
Thus $v$ has no 8-term alternating subsequence. In $v$ immediate repetitions may appear only on 
the transitions $I_iI_{i+1}$. Deleting at most $n-1$ (actually $n-2$ because 
$I_1=\emptyset$) terms from $v$ we obtain a 2-sparse subsequence $w$ on which we can apply 
(\ref{lambda6}). We have
$$
|E|=|v|\le |w|+n-1\le\lambda_6(n)+n-1\ll n\cdot 2^{(1+o(1))\alpha(n)^2}.
$$

On the other hand, let $n\in\N$ and $v$ be the longest 2-sparse sequence over $[n]$ such that 
$v\not\succ ababa$. It uses all $n$ symbols and, by the lower bound in (\ref{lambda3}),
$|v|>cn\alpha(n)$ for an absolute constant $c>0$. Notice that every $i\in[n]$ appears in $v$ at least twice. 
We rename the symbols in $v$ so that for every $1\le i<j\le n$ 
the first appearance of $j$ in $v$ precedes that of $i$; this affects neither the property $v\not\succ ababa$ 
nor the 2-sparseness. By an {\em extremal term of\/} $v$ we mean the first or the last 
appearance of a symbol in $v$. The sequence $v$ has exactly $2n$ extremal terms. We decompose $v$ uniquely 
into intervals $v=I_1I_2\dots I_{2n}$ so that every $I_i$ ends with an extremal term and contains 
no other extremal 
term. Every $I_i$ consists of distinct terms because a repetition 
$\dots b\dots b\dots$ in 
$I_i$ would force a $5$-term alternating subsequence $\dots a\dots b\dots a\dots b\dots a\dots $ in $v$. 
We define a simple (bipartite) graph $G^*=([3n],E)$ by 
$$
\{i,j\}\in E\Longleftrightarrow i\in[n]\;\&\;j\in[n+1,3n]\;\&\;
\mbox{$i$ appears in $I_{j-n}$}.
$$
$G^*$ has $3n$ vertices and $|E|=|v|>cn\alpha(n)$ edges. Suppose that $G^*$ contains the forbidden
ordered subgraph $G_0$ on the vertices $1\le a_1<a_2<\ldots<a_6\le 3n$. By the definition of $G^*$, 
$z=a_1a_2a_1a_2$ is a subsequence of $v$ and its terms appear in $I_{a_3-n},\dots, I_{a_6-n}$, respectively. 
Since $a_2>a_1$, number $a_2$ must appear in $v$ before $z$ starts and therefore $v$ contains a $5$-term 
alternating subsequence but this is forbidden. So $G^*$ does not contain $G_0$ and shows that
$$
g(n)\gg n\alpha(n).
$$ 
This concludes the proof of (\ref{gn}).

\begin{prob}
Narrow the gap $\lambda_3(n)\ll g(n)\ll\lambda_6(n)$ in (\ref{gn}). What is the precise asymptotics of 
$g(n)$?
\end{prob}

Our illustrative example with $g(n)$ shows that the ordered version of a simple graph extremal 
problem may differ dramatically from the unordered one.  
Classical extremal theory of graphs and hypergraphs deals with unordered vertex sets and it produced many 
results of great variety --- see, for example,  
Bollob\'as \cite{boll78, boll95}, Frankl \cite{fran}, 
F\"{u}redi \cite{fure91}, and Tuza \cite{tuza94, tuza96}. However, only little attention has been paid to 
ordered extremal problems. The only systematic studies devoted to this topic known to us 
are F\"uredi and Hajnal \cite{fure_hajn} (ordered bipartite graphs) and Brass, K\'arolyi and Valtr \cite{bras_valt} (cyclically ordered graphs). We think that for several 
reasons ordered extremal problems should be studied and investigated more intensively. First, for their 
intrinsic combinatorial beauty. Second, since they present to us new functions not to be met 
in the classical theory: nearly linear extremal functions, like $n\alpha(n)$ or 
$n\log n$, are characteristic for ordered extremal problems and it seems that they cannot appear without 
ordering of some sort. Third, estimates coming from ordered extremal theory were successfuly applied in 
combinatorial geometry, where often the right key to a problem turns out to be some linear or partial 
ordering, and to obtain more such applications we have to understand more thoroughly combinatorial cores 
of these arguments. 

Before summarizing our results, we return to DS sequences and show that the sequential containment $\prec$ 
can be naturally interpreted in terms of particular hypergraphs, (set) partitions. 
A sequence $u=a_1a_2\dots a_r$ 
over the alphabet $A$ may be viewed as a partition $P$ of $[r]$ where $i$ and $j$ are in the same block of 
$P$ if and only if $a_i=a_j$. Thus blocks of $P$ correspond to the positions of symbols in $u$. 
For example, $u=abaccba$ is the partition $\{\{1,3,7\},\{2,6\},\{4,5\}\}$. 
If 
$u=([r],\sim_u)$ and $v=([s],\sim_v)$ are two sequences given as partitions by equivalence relations, then 
$u\prec v$ if and only if there is an {\em increasing\/} injection $f:[r]\to[s]$ such that 
$x\sim_u y\Leftrightarrow f(x)\sim_v f(y)$ holds for every $x,y\in[r]$.

In this article we introduce and investigate a hypergraph containment that generalizes both the ordered 
subgraph relation and the sequential containment. The containment and its associated extremal 
functions $\ex_e(F,n)$ and $\ex_i(F,n)$ are given in Definitions~\ref{defofcont} and \ref{defoffunc}. 
The function $\ex_e(F,n)$ counts edges in extremal simple hypergraphs $H$ not containing a fixed hypergraph 
$F$ and the function $\ex_i(F,n)$ counts sums of edge cardinalities. In 
Theorem~\ref{exiaexe} we show that for many $F$ one has $\ex_i(F,n)\ll\ex_e(F,n)$. Theorem~\ref{blowups} 
shows that if $F$ is a simple graph, then in some cases good bounds on $\ex_e(F,n)$ can be obtained from bounds 
on the ordered graph extremal function $\mathrm{gex}(F,n)$. We apply Theorem~\ref{blowups} to prove in 
Theorem~\ref{obecfure} that for $G_1=(\{1,3\},\{1,5\},\{2,3\},\{2,4\})$ one has 
$\ex_e(G_1,n)\ll n\cdot(\log n)^2\cdot(\log\log n)^3$ and the same bound for $\ex_i(G_1,n)$;
this generalizes the bound $\gex(G_1,n)\ll n\cdot\log n$ 
of F\"uredi. In another application, Theorem~\ref{klasstro}, we prove that the {\em unordered\/} hypergraph 
extremal function $\ex_e^u(F,n)$ of every forest $F$ is $\ll n$ . In Theorem~\ref{ipomocie} we 
generalize the bound (\ref{mujobec}) to hypergraphs. In Theorem~\ref{starfore} we prove 
that if $F$  is a star forest, then $\ex_e(F,n)$ has an almost linear upper bounds in terms of $\alpha(n)$; 
this generalizes the upper bound in (\ref{gn}). In the concluding section we introduce the notion 
of orderly bipartite forests and pose some problems. 

This article is a revised version of about one half of the material in the technical report \cite{klaz01}. 
We present the other half in \cite{klaz_dalsi}. 

\section{Definitions and bounding weight by size}

A {\em hypergraph\/} $H=(E_i:\ i\in I)$ is a finite list of finite nonempty subsets $E_i$ of 
$\N=\{1,2,\ldots\}$, called {\em edges\/}. $H$ is {\em simple\/} if $E_i\neq E_j$ for every $i,j\in I$, 
$i\neq j$. $H$ is a {\em graph\/} if $|E_i|=2$ for every $i\in I$. $H$ is a {\em partition\/} if 
$E_i\cap E_j=\emptyset$ for every $i,j\in I$, $i\neq j$. The elements 
of $\bigcup H=\bigcup_{i\in I} E_i\subset\N$ are called {\em vertices\/}. Note that our hypergraphs have 
no isolated vertices. The {\em simplification\/} of $H$ is the simple hypergraph obtained from $H$ by 
keeping from each family of mutually equal edges just one edge. 

\begin{defi}\label{defofcont}
Let $H=(E_i:\ i\in I)$ and $H'=(E_i':\ i\in I')$ be two hypergraphs. $H$ {\em contains\/} 
$H'$, in symbols $H\succ H'$, if there exist an {\em increasing\/} injection $F:\bigcup H'\to\bigcup H$ and an 
injection $f: I'\rightarrow I$ such that the implication 
$$
v\in E_i'\Rightarrow F(v)\in E_{f(i)}
$$
holds for every vertex $v\in\bigcup H'$ and every index $i\in I'$. Else we say that $H$ is 
$H'$-{\em free\/} and write $H\not\succ H'$.
\end{defi}

The hypergraph containment $\prec$ clearly extends the sequential containment and the ordered subgraph 
relation. We give an alternative definition of $\prec$. $H=(E_i:\ i\in I)$ and $H'=(E_i':\ i\in I')$ are 
{\em isomorphic\/} if there are an {\em increasing\/} bijection $F:\bigcup H'\to\bigcup H$ and a bijection 
$f: I'\rightarrow I$ such that $F(E_i')=E_{f(i)}$ for every $i\in I'$. $H'$ is a {\em reduction\/} 
of $H$ if $I'\subset I$ and $E_i'\subset E_i$ for every $i\in I'$. Then $H'\prec H$ if and only if $H'$ 
is isomorphic to a reduction of $H$. We call that reduction of $H$ an $H'$-{\em copy\/} in $H$.
For example, if $H'=(\{1\}_1,\{1\}_2)$ ($H'$ is a singleton edge repeated twice) then $H\not\succ H'$
iff $H$ is a partition. Another example: If $H'=(\{1,3\},\{2,4\})$ then $H$ is $H'$-free iff $H$ has no 
four vertices $a<b<c<d$ such that $a$ and $c$ lie in one edge of $H$ and $b$ and $d$ lie in 
another edge.  

The {\em order\/} $v(H)$ of $H=(E_i:\ i\in I)$ is the number of vertices $v(H)=|\bigcup H|$, the 
{\em size\/} $e(H)$ is the number of edges $e(H)=|H|=|I|$, and the {\em weight\/} $i(H)$ is the number of 
incidences between the vertices and the edges $i(H)=\sum_{i\in I}|E_i|$. Trivially, $v(H)\le i(H)$ and 
$e(H)\le i(H)$ for every $H$.

\begin{defi}\label{defoffunc}
Let $F$ be any hypergraph. We associate with $F$ the extremal functions $\ex_e(F),\ex_i(F): \N\to\N$, defined by
\begin{eqnarray*}
\ex_e(F,n)&=&\max\{e(H):\  H\not\succ F\;\&\;\mbox{$H$ is simple}\;\&\;v(H)\le n\}\\
\ex_i(F,n)&=&\max\{i(H):\  H\not\succ F\;\&\;\mbox{$H$ is simple}\;\&\;v(H)\le n\}.
\end{eqnarray*}
\end{defi}

We considered $\ex_e(F,n)$ and $\ex_i(F,n)$ implicitly already in Klazar \cite{klaz00}. Except of this 
article, to our knowledge, this extremal setting is new and was not investigated before. 
Obviously, for every $n\in\N$ and 
$F$, $\ex_e(F,n)\le 2^n-1$ and $\ex_i(F,n)\le n2^{n-1}$ but much better bounds can be usually given.
The {\em reversal\/} of a hypergraph $H=(E_i:\ i\in I)$ with $N=\max(\bigcup H)$ is the hypergraph 
$\overline{H}=(\overline{E_i}:\ i\in I)$ where $\overline{E_i}=\{N-x+1:\ x\in E_i\}$. Thus reversals are 
obtained by reverting the linear ordering of vertices. It is clear that 
$\ex_e(F,n)=\ex_e(\overline{F},n)$ and $\ex_i(F,n)=\ex_i(\overline{F},n)$ for every $F$ and $n$. 

We give a few comments on Definitions~\ref{defofcont} and \ref{defoffunc}. Replacing the implication 
in Definition~\ref{defofcont} with an equivalence, we obtain an induced hypergraph containment that still 
extends the sequential containment. Let $\mathrm{iex}_e(F,n)$ be the corresponding extremal function. For 
$F_k=(\{1\},\{2\},\dots,\{k\})$ and $k\ge 2$ we have $\mathrm{iex}_e(F_k,n)\ge{n\choose k-2}$ 
because the hypergraph $(E:\ E\subset[n]\;\&\;|E|=n-k+2)$ does not contain $F_k$ in the induced sense. 
But in the hypergraph "DS theory" such uncomplicated hypergraphs like $F_k$ should have linear 
(or almost linear) extremal functions. Thus the induced containment is not the right generalization, as far as 
one is interested in "DS theories". 

For graphs, if $H_2\succ H_1$ and $H_2$ is simple then $H_1$ is simple as well. But simple
hypergraphs may contain hypergraphs that are not simple. 
In Definition~\ref{defoffunc} $H$ must be simple because allowing all $H$ would produce usually the value 
$+\infty$ (the simplicity of $H$ may be dropped only for $F=(\{1\}_1,\{1\}_2,\ldots,\{1\}_k)$). On the other 
hand, we allow for the forbidden $F$ any hypergraph: $F$ need not be simple and may have singleton edges. Such a
freedom for $F$ seemed to one referee of \cite{klaz01} "very artificial". This opinion is 
understandable but the author does not share it. Putting aside psychological inertia, there is no reason why 
to restrict $F$ from the outset to be simple or in any other way. On the contrary, doing so we might miss 
some connections and phenomena. Thus we start 
in the definitions with completely arbitrary $F$ and restrict it later only if circumstances require so.

Another perhaps unusual feature of our extremal theory is that in $H$ and $F$ edges of all cardinalities are 
allowed; in extremal theories with forbidden substructures it is more common to have edges of just one 
cardinality. This led naturally to the function $\ex_i(F,n)$ which accounts for edges of all sizes. Trivially, 
$\ex_i(F,n)\ge\ex_e(F,n)$ for every hypergraph $F$ and $n\in\N$. On the other hand, Theorem~\ref{exiaexe} shows 
that for many $F$ one has $\ex_i(F,n)\ll\ex_e(F,n)$. In Definition~\ref{defoffunc} we take all $H$ with 
$v(H)\le n$ that the extremal functions be automatically nondecreasing. Replacing $v(H)\le n$ with $v(H)=n$ 
would give more information on the functions but also would bring the complication that then extremal 
functions are not always nondecreasing. It happens for $F=(\{1\},\{2\},\dots,\{k\})$ and we analyze this 
phenomenon in \cite{klaz01,klaz_dalsi}. 

\begin{veta}\label{exiaexe}
Suppose that the hypergraph $F$ has no two separated edges, which means that $E_1<E_2$ holds for no two 
edges of $F$. Let 
$p=v(F)$ and $q=e(F)>1$. Then for every $n\in\N$, 
$$
\ex_i(F,n)\le (2p-1)(q-1)\cdot\ex_e(F,n).
$$
\end{veta}
\duk
Suppose that $H$ attains the value $\ex_i(F,n)$. We transform $H$ in a new hypergraph $H'$ by keeping all 
edges with less than $p$ vertices and replacing every edge $E=\{v_1,v_2,\ldots,v_s\}$ of $H$  
with $s\ge p$, where $v_1<v_2<\cdots<v_s$, by $t=\lfloor |E|/p\rfloor$ new $p$-element edges 
$\{v_1,\ldots,v_p\}$, $\{v_{p+1},\ldots,v_{2p}\},\ldots,$ $\{v_{(t-1)p+1},\ldots,v_{tp}\}$. $H'$ may not be 
simple and we set $H''$ to be the simplification of $H'$. Two observations: (i) no edge of $H'$ repeats $q$ 
or more times and (ii) $H''$ is $F$-free. If (i) were false, there would be $q$ distinct edges 
$E_1,\ldots,E_q$ in $H$ such that $|\bigcap_{i=1}^q E_i|\ge p$. But this implies the contradiction 
$F\prec H$. As for (ii), any $F$-copy in $H''$ may use from every $E\in H$ at most one new edge 
$E''\subset E$ (the new edges born from $E$ are mutually separated) and so it is an $F$-copy in $H$ as well. 
The observations and the definitions of $H'$ and $H''$ imply
\begin{eqnarray*}
\ex_i(F,n)=i(H)&\le& {(2p-1)\cdot i(H')\over p}\le {(2p-1)(q-1)\cdot i(H'')\over p}\\
&\le& (2p-1)(q-1)\cdot e(H'')\\
&\le& (2p-1)(q-1)\cdot\ex_e(F,n).
\end{eqnarray*}
The last innocently looking inequality follows from the fact that $\ex_e(F,n)$ is nondecreasing by definition.
\kduk

However, $\ex_i(F,n)\ll\ex_e(F,n)$ does not hold for $F_k=(\{1\},\{2\},\dots,\{k\})$ and $k\ge 2$: 
$\ex_e(F_k,n)=2^{k-1}-1$ for $n\ge k-1$ and $\ex_i(F_k,n)=(k-1)n-(k-2)$ for $n>\max(k,2^{k-2})$
(\cite{klaz01,klaz_dalsi}). Note that for $F=(\{1\})$ both extremal functions are undefined. $F_k$ is highly 
symmetric and the ordering of vertices is irrelevant for the containment $H\succ F_k$. 

\begin{prob}
Prove that if $F$ is not isomorphic to $(\{1\},\{2\},\dots,\{k\})$, $k\ge 1$, then $\ex_i(F,n)\ll\ex_e(F,n)$. 
\end{prob}

\section{Bounding hypergraphs by means of graphs}
For a family of simple graphs $R$ and $n\in\N$ we define
\begin{eqnarray*}
\gex(R,n)&=&\max\{e(G):\ G\not\succ G' \mbox{ for all } G'\in R\;\&\;\mbox{$G$ is a simple 
graph}\\
&&\;\&\;v(G)\le n\}
\end{eqnarray*}
and for one simple graph $G$ we write $\gex(G,n)$ instead of $\gex(\{G\},n)$. F\"uredi proved in \cite{fure90}, 
see also \cite{fure_hajn}, that for
\begin{equation}\label{GF}
G_1=(\{1,3\},\{1,5\},\{2,3\},\{2,4\})=\ 
\obr{\put(1,1){\bod}\put(8,1){\bod}\put(15,1){\bod} 
\put(22,1){\bod}\put(29,1){\bod}
\put(8,1){\line(1,0){7}}
\put(8,1){\oval(14,10)[t]}\put(15,1){\oval(14,7)[t]}
\put(15,1){\oval(28,16)[t]}
}
\end{equation}
one has 
\begin{equation}\label{furebound}
n\log n\ll\gex(G_1,n)\ll n\log n. 
\end{equation}
(In \cite{fure90} and \cite{fure_hajn} the investigated objects are 0-1 matrices, which are ordered bipartite 
graphs, but in the case of $G_1$ the bounds are easily extended to all ordered graphs.) Attempts to generalize 
the upper bound in (\ref{furebound}) to hypergraphs motivated the next theorem. 

For $k\in\N$ we say that a simple graph $G'$ is a $k$-{\em blow-up\/} of a simple graph $G$ if for 
every edge coloring 
$\chi: G'\to\N$ that uses every color at most $k$ times there is a $G$-copy in $G'$ with totally different
colors, that is, $\chi$ is injective on the $G$-copy. For $k\in\N$ and a simple graph $G$ we write $B(k,G)$ to 
denote the set of all $k$-blow-ups of $G$.

\begin{veta}\label{blowups}
Let $F$ be a simple graph with $p=v(F)$ and $q=e(F)>1$ and let $B\subset B({p\choose 2},F)$.
If $f: \N\to\N$ is a nondecreasing function such that 
$$
\mathrm{gex}(B,n)<n\cdot f(n)
$$
holds for every $n\in\N$, then
\begin{equation}\label{rekner}
\ex_e(F,n)< q\cdot \mathrm{gex}(F,n)\cdot \ex_e(F,2f(n)+1)
\end{equation}
holds for every $n\in\N$, $n\ge 3$.
\end{veta}
\duk
Suppose that the simple hypergraph $H$ attains the value $\ex_e(F,n)$ and $\bigcup H=[m]$, $m\le n$. We put in 
$H'$ every edge of $H$ with more than $1$ and less than $p$ vertices and for every $E\in H$ with $|E|\ge p$ 
we put in $H'$ an arbitrary subset $E'\subset E$, $|E'|=p$. So $2\le |E|\le p$ for every $E\in H'$ and no edge 
of $H'$ repeats more then $q-1$ times, for else we would have $H\succ F$. Let $H''$ be the simplification 
of $H'$. We have
$$
e(H)\le n+(q-1)e(H'').
$$ 
Let $G$ be the simple graph consisting 
of all edges $E^*$ such that $E^*\subset E$ for some $E\in H''$. Observe that if $F'\in B$ and $F'\prec G$, 
then $F\prec H''$ and thus $F\prec H$. (For the edges $E^*\in G$ forming an $F'$-copy consider the coloring 
$\chi(E^*)=E$ where $E\in H''$ is such that 
$E^*\subset E$. Every color is used at most ${p\choose 2}$ times and therefore, since $F'$ is a 
${p\choose 2}$-blow-up of $F$, we have an $F$-copy in $G$ for which the correspondence $E^*\mapsto E$ is 
injective.) 
Hence $F'\prec G$ for no $F'\in B$. Let $v(G)=n'$; $n'\le n$. We have 
$$
e(G)\le \mathrm{gex}(B,n')<n'\cdot f(n').
$$

There exists a vertex $v_0\in\bigcup G$ such that 
$$
d=\deg_{G}(v_0)<2f(n')\le 2f(n).
$$
Fix an arbitrary edge $E_0^*$, $v_0\in E_0^*\in G$. Let $X\subset[n]$ be the union of all edges $E\in H''$ 
satisfying $E_0^*\subset E$ and $m$ be the number of such edges in $H''$. We have the inequalities 
$$
m\le\ex_e(F,|X|)\ \mbox{ and }\ |X|\le d+1. 
$$
Thus
$$
m\le \ex_e(F,|X|)\le \ex_e(F,d+1)\le\ex_e(F,2f(n)+1). 
$$
We see that the two-element set $E_0^*$ is contained in at least 1 but at most $\ex_e(F,2f(n)+1)$ edges of 
$H''$. Deleting those edges we obtain a subhypergraph $H_1''$ of $H''$ on which the same argument can be 
applied. That is, a two-element set $E_1^*$ exists such that $E_1^*\subset E$ for at least 1 but at most 
$\ex_e(F,2f(n)+1)$ edges $E\in H_1''$ and clearly $E_1^*\neq E_0^*$. Continuing this way until the whole 
$H''$ is exhausted, we define a mapping
$$
M: H''\to\{E^*:\ E^*\subset [n], |E^*|=2\}
$$
such that
$$
M(E)\subset E\ \mbox{ and }\ 
|M^{-1}(E^*)|\le\ex_e(F,2f(n)+1)
$$
holds for every $E\in H''$ and $E^*\subset [n], |E^*|=2$. Let $G'$ be the simple graph $G'=M(H'')$
and $v(G')=n'$; $n'\le n$.

The containment $F\prec G'$ implies, by the definition of $G'$, that 
$F\prec H''$ and hence $F\prec H$, which is not allowed. Thus 
$$
e(G')\le \mathrm{gex}(F,n')\le \mathrm{gex}(F,n).
$$ 
Putting it all together, we obtain (since $\mathrm{gex}(F,n)\ge n-1$ if $q>1$)
\begin{eqnarray*}
\ex_e(F,n)=e(H)&\le & n+(q-1)\cdot e(H'')\\
 &\le& n+(q-1)\cdot\ex_e(F,2f(n)+1)\cdot e(G')\\
&<& q\cdot \ex_e(F,2f(n)+1)\cdot \mathrm{gex}(F,n)
\end{eqnarray*}
for every $n\ge 3$.
\kduk

\noindent
We give three applications of this theorem. The first one is the promised generalization of the upper bound 
in (\ref{furebound}).  

We say that a simple graph $G$ is a $k$-{\em multiple\/} of $G_1$, where $k\in\N$ and $G_1$ is defined 
in (\ref{GF}), if $G$ has this structure: $\bigcup G=A\cup\{v\}\cup B\cup C$ with $A<v<B<C$, $|A|=k$, the 
vertex $v$ has degree $k$ and is connected to every vertex in $A$, every vertex in $A$ has degree $2k+1$ and 
is besides $v$ connected to $k$ vertices in $B$ and to $k$ vertices in $C$, and $G$ has no other edges. The 
edges incident with $v$ are called {\em backward edges\/} and the edges incident with vertices in $B\cup C$ 
are called {\em forward edges\/}. We denote the set of all $k$-multiples of $G_1$ by $M(k)$. 

\begin{lema}\label{blowupfure}
The sets of graphs $M(k)$, $k\in\N$, have the following properties.  
\begin{enumerate}
\item For every $k$, $M(3k+1)\subset B(k,G_1)$. In particular, $M(31)\subset B({5\choose 2},G_1)$.
\item For every $k$, $\gex(M(k),n)\ll n\log n$.
\end{enumerate}
\end{lema}
\duk
1. Let $G$ be a $3k+1$-multiple of $G_1$ and $\chi:G\to\N$ be an edge coloring using each color at most 
$k$ times. We select in $G$ two backward edges $E_1=\{i,v\}$ and $E_2=\{j,v\}$, $i<j<v$, with different 
colors. It follows that we can select in $G$ two forward edges $E_3=\{i,l\}$ and $E_4=\{j,l'\}$ 
such that $v<l'<l$ and the colors $\chi(E_1),\dots,\chi(E_4)$ are distinct. Edges $E_1,\dots,E_4$ form a 
$G_1$-copy on which $\chi$ is injective. Thus $G\in B(k,G_1)$. 

2. Let $n\ge 2$ and $G$ be any simple graph such that $\bigcup G=[n]$ and $G\not\succ F$ for every $F\in M(k)$. 
Let $J_i=\{E\in G:\ \min E=i\}$ for $i\in[n]$. In every $J_i$ it is possible to mark 
$\lfloor|J_i|/k\rfloor>|J_i|/k-1$ edges so that every marked edge is in $J_i$ immediately 
followed by $k-1$ unmarked ones (we order the edges in $J_i$ by their endpoints). The graph $G'$ formed by the 
marked edges satisfies 
$$
e(G')> e(G)/k-n.
$$ 
Also, for every edge $\{i,j\}\in G'$, $i<j$, there are 
at least $k-1$ edges $\{i,l\}\in G$ with $l>j$, and for every two edges $\{i,j\},\{i,j'\}\in G'$, $i<j<j'$,
there are at least $k-1$ edges $\{i,l\}\in G$ with $j<l<j'$.
Now we proceed as in F\"uredi \cite{fure90}. We say that $\{i,j\}\in G'$, $i<j$, has {\em type\/} $(j,m)$, 
where $m\ge 0$ is an integer, if 
there are two edges $\{i,l\}$ and $\{i,l'\}$ in $G'$ such that $j<l<l'$ and $l-j\le 2^m<l'-j$. 
Consider the partition 
$$
G'=G^*\cup G^{**}
$$ 
where $G^*$ is formed by edges with at least one type and $G^{**}$ by 
edges without type. It follows from the definition of type and of $G'$ that if $k$ edges of $G^*$ have 
the same type, then $F\prec G$ for some $F\in M(k)$ which is forbidden. Thus any type is shared by at 
most $k-1$ edges. Since the number of types is less than $n(1+\log_2 n)$, we have
$$
e(G^*)<(k-1)n+(k-1)n\log_2n.
$$ 
To bound $e(G^{**})$, we fix a vertex $i\in[n]$ and consider the endpoints $i<j_0<j_1<\cdots<j_{t-1}\le n$ 
of all edges $E\in G'$ which have no type and $\min E=i$. Let $d_r=j_r-j_{r-1}$ for $1\le r\le t-1$ and 
$D=d_1+\cdots+d_{t-1}=j_{t-1}-j_0$. If $d_1\le D/2$, then $d_1\le 2^m<D$ for 
some integer $m\ge 0$ and the edge $\{i,j_0\}$ would have type $(j_0,m)$ because of the edges $\{i,j_1\}$ and 
$\{i,j_{t-1}\}$. Thus $d_1>D/2$ and $D-d_1<D/2$. By the same  argument applied to $\{i,j_1\}$, 
$d_2>(D-d_1)/2$ and thus $D-d_1-d_2<D/4$. 
In general, $1\le D-d_1-\cdots-d_r<D/2^r$ for $1\le r\le t-2$. Thus 
$t\le \log_2 D+2<2+\log_2n$. Summing these inequalities for all $i\in[n]$, we have 
$$
e(G^{**})<2n+n\log_2n.
$$ 
Alltogether we have 
$$
e(G)<kn+k(e(G^*)+e(G^{**}))<(k^2+k)n+k^2n\log_2 n. 
$$
We conclude that $\gex(M(k),n)\ll n\log n$ and the constant in $\ll$ depends quadratically on $k$.
\kduk

\begin{veta}\label{obecfure}
Let $G_1$ be the simple graph given in (\ref{GF}). We have the following bounds. 
\begin{enumerate}
\item $n\cdot\log n\ll \ex_e(G_1,n)\ll n\cdot(\log n)^2\cdot(\log\log n)^3$.
\item $n\cdot\log n\ll \ex_i(G_1,n)\ll n\cdot(\log n)^2\cdot(\log\log n)^3$.
\end{enumerate}
\end{veta}
\duk
1. The lower bound follows from the lower bound in (\ref{furebound}). To prove the upper bound, we use 
Theorem~\ref{blowups}. By 2 of Lemma~\ref{blowupfure}, we have 
$\gex(M(31),n)\ll n\log n$. Also, $\gex(G_1,n)\ll n\log n$ (by the upper bound in (\ref{furebound}) or
by $\gex(G_1,n)\le\gex(B,n)$). By 1 of Lemma~\ref{blowupfure}, we can apply Theorem~\ref{blowups} with 
$B=M(31)$. Starting with the trivial bound $\ex_e(G_1,n)<2^n$, (\ref{rekner}) with 
$f(n)\ll\log n$ gives 
$$
\ex_e(G_1,n)\ll n^c
$$
where $c>0$ is a constant. Feeding this bound back to (\ref{rekner}), we get
$$
\ex_e(G_1,n)\ll n\cdot(\log n)^{c+1}.
$$
Two more iterations of (\ref{rekner}) give 
$$
\ex_e(G_1,n)\ll n\cdot(\log n)^2\cdot(\log\log n)^{c+1}
$$
and
$$
\ex_e(G_1,n)\ll n\cdot(\log n)^2\cdot(\log\log n)^2\cdot(\log\log\log n)^{c+1}
$$
which is slightly better than the stated bound. 

2. The lower bound follows from $\ex_i(G_1,n)\ge\ex_e(G_1,n)$. The upper bound follows from the upper 
bound in 1 by Theorem~\ref{exiaexe}. 
\kduk

\begin{prob}
What is the exact asymptotics of $\ex_e(G_1,n)$?
\end{prob}

The second application of Theorem~\ref{blowups} concerns unordered extremal functions $\ex^u_e(F,n)$ and 
$\gex^u(G,n)$. They are 
defined as $\ex_e(F,n)$ and $\gex(G,n)$ except that in the containment the injection need not be 
increasing. So $\gex^u(G,n)$ is the classical graph extremal function. It is well known, see for example 
Bollob\'as \cite[Exercise 24 in IV.7]{boll98}, that $\gex^u(F,n)\le (e(F)-1)\cdot n$ for every forest $F$. 
We extend the linear bound to unordered hypergraphs. Theorem~\ref{blowups} holds also in the 
unordered case because the proof is independent of ordering. Ordering is crucial only for obtaining 
linear or almost linear bounds on $\gex(F,n)$ and $\gex(B,n)$ because the inequality (\ref{rekner}) is useless 
if $f(n)$ is not almost constant. The proof of Theorem~\ref{blowups} shows also that if $F$ is a forest 
and all members of $B$ are forests (which is not the case for $B=M(k)$) then ${p\choose 2}$ can be 
replaced by $p-1$ (because for $|E|=p$ every $p$ 
two-element edges $E^*\subset E$ contain a cycle but no $F'\in B$ has a cycle).

\begin{veta}\label{klasstro}
Let $F$ be a forest. Its unordered hypergraph extremal function satisfies
$$
\ex_e^u(F,n)\ll n.
$$
\end{veta}
\duk
Let $v(F)=p$ and $e(F)=q>1$ (case $q=1$ is trivial). It is not hard to find a forest 
$F'$ with large $e(F')=Q$ --- $Q\le (pq^2)^{q+1}$ suffices --- 
that is a $p-1$-blow-up of $F$. We set $B=\{F'\}$ 
and use (\ref{rekner}) with the bounds $\gex^u(F,n)\le (q-1)n$, $f(n)=Q-1$ (since 
$\gex^u(B,n)=\gex^u(F',n)\le (Q-1)n$), and $\ex_e^u(F,n)<2^n$ (trivial):
$$
\ex_e^u(F,n)<q\cdot (q-1)n\cdot 2^{2Q-1}=
{\textstyle {q\choose 2}}4^Q\cdot n.
$$
\kduk

\noindent
One can prove the bound $\ex^u_e(F,n)\ll n$ also directly, without Theorem~\ref{blowups}, by adapting the proof 
of $\gex^u(F,n)\ll n$ to hypergraphs.
The third application of Theorem~\ref{blowups} follows in the next section. 

\section{Partitions and star forests}

The bound (\ref{mujobec}) tells us that if $F$ is any fixed partition 
with $k$ blocks and $H$ is a $k$-sparse partition with $H\not\succ F$, then $v(H)=i(H)$ has an almost linear 
upper bound in terms of $e(H)$. The following theorem bounds $i(H)$ almost linearly in terms of 
$e(H)$ in the wider class of (not necessarily simple) hypergraphs $H$. 
The proof is based on (\ref{muj}).

\begin{veta}\label{ipomocie}
Let $F$ be a partition with $p=v(F)$ and $q=e(F)>1$ and $H$ be a $F$-free hypergraph, not 
necessarily simple. Then
\begin{equation}\label{iae}
i(H)<(q-1)\cdot v(H)+e(H)\cdot\beta(q,2p,e(H))
\end{equation}
where $\beta(k,l,n)$ is the almost constant function defined in (\ref{beta}). 
\end{veta}
\duk
Let $\bigcup H=[n]$ and the edges of $H$ be $E_1,E_2,\ldots,E_e$ where $e=e(H)$. 
We set, for $1\le i\le n$, $S_i=\{j\in[e]:\ i\in E_j\}$ and consider the sequence 
$$
v=I_1I_2\ldots I_n
$$ 
where $I_i$ is an arbitrary ordering of $S_i$. Clearly, 
$v$ is over $[e]$ and $|v|=i(H)$. The sequence $v$ may not be $q$-sparse, because of the transitions 
$I_iI_{i+1}$, but it is easy to delete at most $q-1$ terms from the beginning of 
every $I_i$, $i>1$, so that the resulting subsequence $w$ is $q$-sparse. Thus $|w|\ge |v|-(q-1)(n-1)$. It 
follows that 
if $w$ contains $u(q,2p)$, where $u(k,l)$ is defined in (\ref{ukl}), then $H$ contains $F$ but
this is forbidden. (Note that the 
subsequence $aab$ in $v$ forces the first $a$ and the $b$ to appear in two distinct segments 
$I_i$ and thus it gives incidences of $E_a$ and $E_b$ with two distinct vertices.) Hence $w\not\succ u(q,2p)$ 
and we can apply (\ref{muj}):
$$
i(H)=|v|<(q-1)n+|w|\le (q-1)n+e\cdot\beta(q,2p,e).
$$
\kduk

\noindent
We show that for the partition
$$
F=H_2=(\{1,3,5\},\{2,4\})=\ 
\obr{\put(3,1){\bod}\put(9,1){\bod}\put(15,1){\bod} 
\put(21,1){\bod}\put(27,1){\bod}
\put(15,1){\oval(30,14)[t]}
\put(3,1){\oval(6,10)[b]}
\put(27,1){\oval(6,10)[b]}
\put(9,1){\oval(6,8)[t]}
\put(15,1){\oval(6,10)[b]}
\put(21,1){\oval(6,8)[t]}
\put(15,1){\oval(12,5)[b]}
}
$$

\bigskip\noindent
the factor at $e(H)$ in (\ref{iae}) must be $\gg\alpha(e(H))$. We proceed as in the proof of 
$g(n)=\gex(G_0,n)\gg n\alpha(n)$ in (\ref{gn}) and take a 2-sparse sequence $v$ over $[n]$ such that 
$v\not\succ 12121$, $|v|\gg n\alpha(n)$, and $v=I_1I_2\ldots I_{2n}$ where every interval $I_i$ 
consists of distinct terms. We define the hypergraph 
$$
H=(E_i:\ i\in[n])\ \mbox{ with }\ E_i=\{j\in[2n]:\ i \mbox{ appears in }\  I_j\}.
$$
We have $i(H)=|v|\gg n\alpha(n)$, $\bigcup H=[2n]$, $v(H)=2n$, and $e(H)=n$. It is clear that 
$H\not\succ H_2$ because $v\not\succ 12121$.

Taking in Theorem~\ref{ipomocie} $H$ to be simple and with the maximum weight, 
we obtain as a corollary that if $F$ is a partition, $p=v(F)$, and $q=e(F)>1$, then 
$$
\ex_i(F,n)< (q-1)n+\ex_e(F,n)\cdot\beta(q,2p,\ex_e(F,n)). 
$$
But Theorem~\ref{exiaexe}, when it applies, gives better bounds.

Our last theorem generalizes in two ways the upper bound in (\ref{gn}). We consider a class of forbidden 
forests that contains $G_0$ as a member and we extend by means of Theorems~\ref{blowups} and \ref{exiaexe} 
the almost linear upper bound to hypergraphs. The class consists of {\em star forests\/} which are forests $F$ 
with this structure: 
$\bigcup F=A\cup B$ for some sets $A<B$ such that every vertex in $B$ has degree $1$ and every edge of $F$ 
connects
$A$ and $B$. Thus $F$ is a star forest iff every component of $F$ is a star and every central vertex of a star 
is smaller than every leaf.

\begin{veta}\label{starfore}
Let $F$ be a star forest with $r>1$ components, $p$ vertices, and $q=p-r$ edges. Let $t=(p-1)(q-1)+1$ and 
$\beta(k,l,n)$ be the almost constant function defined in (\ref{beta}). We have 
the following bounds. 
\begin{enumerate}
\item $\gex(F,n)<(r-1)n+n\cdot\beta(r,2q,n)$.
\item $\ex_e(F,n)\ll n\cdot\beta(r,2tq,n)^3$.
\item $\ex_i(F,n)\ll n\cdot\beta(r,2tq,n)^3$.
\end{enumerate}
\end{veta}
\duk
1. We mark the centers of the stars in $F$ with $1,2,\dots,r$ according 
to their order and give the leaves of any star the mark of its center. The marks on all leaves then 
form a  
sequence $u$ over $[r]$ of length $p-r$. Now let $G$ be any simple graph with $\bigcup G=[n]$ and 
$G\not\succ F$. We consider the sequence 
$$
v=I_1I_2\ldots I_n
$$
where $I_j$ is any ordering of the set $\{i\in[n]:\ \{i,j\}\in G, i<j\}$. As in the previous proof, we select 
an $r$-sparse subsequence $w$ of $v$ with length 
$|w|\ge |v|-(r-1)(n-1)$. Suppose that $w\succ u(r,2(p-r))$ where $u(k,l)$ is defined in (\ref{ukl}). Then 
$w$ has a (not necessarily consecutive) subsequence $z$ of the form
$$
a_1a_2\dots a_ra_1a_2\dots a_r\dots a_1a_2\dots a_r
$$
with $2(p-r)$ segments $a_1a_2\dots a_r$. We have $a_{i_1}<a_{i_2}<\dots<a_{i_r}$ for a permutation 
$i_1,i_2,\dots,i_r$ of $[r]$. We label every term $a_{i_j}$ in $z$ with $j$. Clearly, if we select 
one term from the 2-nd, 4-th, $\dots$, $2(p-r)$-th segment of $z$ so that the labels on the selected terms 
coincide with the sequence $u$, the selected terms lie in $p-r$ distinct intervals 
$I_{j_1},\dots,I_{j_{p-r}}$, $j_1<\dots<j_{p-r}$. Since the selected terms are preceded by one segment 
$a_1a_2\dots a_r$, we have $a_{i_r}<j_1$. The edges between $a_1,\dots,a_r$ and $j_1,\dots,j_{p-r}$ which 
corespond to the selected terms form an $F$-copy in $G$ but $G\succ F$ is forbidden. Therefore 
$w\not\succ u(r,2(p-r))$ and we can apply (\ref{muj}):
$$
e(G)=|v|\le (r-1)n+|w|<(r-1)n+n\cdot\beta(r,2(p-r),n).
$$

2. Suppose that $F$ has the vertex set $[p]$ (so that $[r]$ are the centers of the stars and $[r+1,p]$ are the 
leaves). For $k\in\N$ we denote $F(k)$ the star forest with the vertex set $[r+(p-r)k]$ in which $[r]$ 
are again the centers 
of stars and for $i=1,2,\dots,p-r$ the vertices in $[r+(i-1)k+1,r+ik]$ are joined to the same vertex 
in $[r]$ as $r+i$ is joined in $F$. 
It is easy to see that $F(t)=F((p-1)(q-1)+1)$ is a $p-1$-blow-up of $F$. Also, $e(F(k))=kq$. 
We set $B=\{F(t)\}$
and use (\ref{rekner}) with the bounds $\gex(F,n)\ll n\cdot\beta(r,2q,n)=n\cdot\beta'$ (bound 1 for 
$F$), $f(n)=c\beta(r,2tq,n)=c\beta$ for a constant $c>0$ (bound 1 for $F(t)$), and $\ex_e(F,n)<2^n$ 
(trivial):
$$
\ex_e(F,n)\ll n\cdot\beta'\cdot 2^{2c\beta+1}<n\cdot 2^{2(c+1)\beta}. 
$$ 
The second application of (\ref{rekner}) gives
$$
\ex_e(F,n)\ll n\cdot\beta'\cdot\beta\cdot 2^{2(c+1)\cdot\beta(r,2tq,2c\beta+1)}\ll n\cdot\beta^3
$$
because $\beta'\le\beta$ and
$$
\beta(r,2tq,x)\ll\log\log x
$$
(this is true with any number of logarithms). 

3. This follows from 2 by Theorem~\ref{exiaexe}.
\kduk

\noindent
The lower bound in (\ref{gn}) shows that in general the factor at $n$ in 1, 2, and 3 of the previous theorem 
cannot be replaced with a constant and may be as big as $\gg\alpha(n)$. The proved bounds hold also for the 
reversals of star forests. 

\section{Concluding remarks}

It is reasonable to call a function $f:\N\to\R$ {\em nearly linear\/} if 
$n^{1-\varepsilon}\ll f(n)\ll n^{1+\varepsilon}$ holds for every $\varepsilon>0$. We identify a candidate for 
the class of all hypergraphs $F$ with nearly linear $\ex_e(F,n)$. If $F$ is isomorphic to the hypergraph
$(\{1\},\{2\},\dots,\{k\})$, then $\ex_e(F,n)$ is eventually constant (\cite{klaz_dalsi}) and thus is not nearly 
linear. For other hypergraphs we have $\ex_e(F,n)\ge n$ because $F\not\prec(\{1\},\{2\},\dots,\{n\})$. An 
{\em orderly bipartite forest\/} is a simple graph $F$ such that $F$ has no cycle and 
$\min E<\max E'$ holds for every two edges of $F$. In other words, $F$ is a forest and there is a partition 
$\bigcup F=A\cup B$ such that $A<B$ and every edge of $F$ connects $A$ and $B$. We denote the class of all 
orderly bipartite forests by OBF. We say that $F$ is an {\em orderly bipartite forest with 
singletons\/} if $F=F_1\cup F_2$ where $F_1\in\mathrm{OBF}$ and $F_2$ is a hypergraph consisting of possibly 
repeating singleton edges. For example, $F$ may be
$$
F=(\{8\},\{6\}_1,\{6\}_2,\{2\},\{1,6\},\{3,6\},\{4,5\},\{4,7\}).
$$
The class OBF subsumes star forests and their reversals. $G_1$ defined in (\ref{GF}) belongs to OBF but is 
neither a star forest nor a reversed star forest. 

\begin{lema}
If the hypergraph $F$ is not an orderly bipartite forest with singletons, then there is a constant 
$\gamma>1$ such that
$$
\ex_e(F,n)\gg n^{\gamma}
$$
and hence $\ex_e(F,n)$ is not nearly linear.
\end{lema}
\duk
If $F$ is not an orderly bipartite forest with singletons, then $F$ has (i) an edge with more than two elements 
or (ii) two separated two-element edges or (iii) a two-path isomorphic to $(\{1,2\},\{2,3\})$ or (iv) a repeated 
two-element edge or (v) an even cycle of two-element edges (odd cycles 
are subsumed in (iii)). In the cases (i)--(iv) we have $\ex_e(F,n)\gg n^2$ because the complete bipartite graph 
with parts $[\lfloor n/2\rfloor]$ and $[\lfloor n/2\rfloor+1,n]$ does not contain $F$. As for the case (v), 
an application of the probabilistic method (Erd\H os \cite{erdo59}) provides 
an unordered graph that has $n$ vertices, $\gg n^{1+1/k}$ edges, and no even cycle of length $k$. Thus, 
in the case (v), $\ex_e(F,n)\gg n^{1+1/k}$ for some $k\in\N$. 
\kduk

\noindent
We conjecture that $\ex_e(F,n)$ is nearly linear if and only if $F$ is an orderly bipartite forest with 
singletons that is not isomorphic to $(\{1\},\{2\},\dots,\{k\})$. Since every orderly bipartite forest with 
singletons is contained in some orderly bipartite forest, it suffices to consider only orderly bipartite 
forests. 

\begin{prob}
Prove (or disprove) that for every orderly bipartite forest $F$ we have 
$$
\ex_e(F,n)\ll n(\log n)^c
$$
for some constant $c>0$.
\end{prob}

\noindent
It is not difficult to find for every $F\in\mathrm{OBF}$ and $k\in\N$ an $F'\in\mathrm{OBF}$ that is a 
$k$-blow-up of $F$. Thus the 
previous bound would follow by Theorem~\ref{blowups} from the graph bound $\gex(F,n)\ll n(\log n)^c$.  

It is natural to consider two subclasses $\mathrm{OBF}^l\subset\mathrm{OBF}^{\alpha}\subset\mathrm{OBF}$ 
where $\mathrm{OBF}^l$ consists of all $F\in\mathrm{OBF}$ with $\ex_e(F,n)\ll n$ and $\mathrm{OBF}^{\alpha}$
consists of all $F\in\mathrm{OBF}$ with $\ex_e(F,n)\ll n\cdot f(\alpha(n))$ for a primitive recursive 
function $f(n)$. Both inclusions are strict as witnessed by $G_0$ and $G_1$ (defined in (\ref{G0}) and 
(\ref{GF})). In this article we ignored the class $\mathrm{OBF}^l$ completely and showed that 
$\mathrm{OBF}^{\alpha}$ contains all star forests (and their reversals). It would much interesting to 
learn more about $\mathrm{OBF}^l$ and $\mathrm{OBF}^{\alpha}$. Does the latter class consist only of 
star forests and their reversals?

\end{document}